\documentclass[times,review,3p]{elsarticle}
\usepackage{graphicx}
\usepackage{amssymb}
\usepackage{amsmath}
\usepackage{hyperref}
\usepackage{lineno}
\usepackage{subfigure}
\usepackage{tikz}
\usepackage{catchfilebetweentags}
\usepackage{booktabs}
\usepackage{multirow}

\modulolinenumbers[1]

\journal{Computers \& Fluids}

\begin{document}

\begin{frontmatter}

\title{A comparative study of split advection algorithms on Moment-of-Fluid (MOF) method  for incompressible flow}
    \author[1]{Zhouteng Ye}
    \ead{yzt9zju@gmail.com}

    \address[1]{
        Ocean College, Zhejiang University, Zhoushan 316021, Zhejiang, People’s Republic of China}

\begin{abstract}
    The moment-of-fluid method (MOF) is known as an extension of
    the volume-of-fluid method with piecewise linear interface construction (VOF-PLIC).
    In this study,
    several directional splitting advection algorithms are extended from the VOF-PLIC
    method to the MOF method.
    Those methods, along with some existing directional splitting algorithms,
    are presented in detail, especially on the geometrical and non-geometrical nature
    of the advection algorithm.
    Besides, we also proposed a simple and efficient analytic form to calculate
    the volume and centroid of the intersection polyhedron in 3D.
    Several numerical tests, including both 2D and 3D tests,
    are carried out to investigate mass conservation and
    geometrical error.
    Numerical results suggest that the mixed EI and LE scheme (EILE2D) \citep{aulisa_geometrical_2003}
    has the best overall performance in 2D and Weymouth and Yue's scheme (WY)
    has the best overall performance in 3D.

\end{abstract}


\end{frontmatter}

\section{Introduction}
\label{sec:intro}
The volume of fluid method (VOF) has become
one of the most popular techniques in tracking the interface between different materials.
In the part 45 years,
a variety of VOF methods has been developed \citep{ehlers_slic_1976,hirt_volume_1981,youngs_time-dependent_1982,zhang_new_2008,zhang_new_2008,maric_unstructured_2020}
Among those VOF methods,
the volume-of-fluid method with piece-wise line interface construction (VOF-PLIC) is one
of the most widely used in tracking the interface within the Eulerian framework.
\citet{dyadechko_moment--fluid_2005} developed Moment-of-Fluid method (MOF) as an
extension of VOF-PLIC method.
With the centroid of the fluid as an additional constraint,
the MOF method resolves the interface with a smaller minimum scale than the VOF-PLIC algorithm
and has been extended from Cartesian grid to multiple frameworks such as adaptive mesh refinement(AMR)\cite{ahn_multi-material_2007, jemison_coupled_2013,liu_moment--fluid_2020},
arbitrary Lagrangian-Eulerian (ALE) \cite{galera_2d_2011,breil_multi-material_2013}.

The MOF method contains two parts:
(1) reconstruction algorithm
and, (2) advection algorithm.
The reconstruction algorithm in MOF method is to use the known reference centroid
and volume fraction to find the optimized linear cut-off in a cell.
Unlike VOF-PLIC method,
the MOF reconstruction does not need the information from the neighbouring grids,
which provides a better resolution when there is insufficient from the neighbouring grid.
Since the \citet{dyadechko_moment--fluid_2005} proposed the first MOF reconstruction algorithm,
several improvements has been done
\citep{ahn_multi-material_2007,jemison_coupled_2013,chen_improved_2016,lemoine_moment--fluid_2017,milcent_moment--fluid_2020}.
In 2D Cartesian grid, the analytic solution to the reconstruction is proposed by
\citet{lemoine_moment--fluid_2017}.
Unfortunately, the 2D analytic solution cannot be extended to 3D.
In 3D,
\citet{milcent_moment--fluid_2020} proposed an analytic form of
the gradient when calculating the normal vector in MOF reconstruction.
Unlike the full analytic solution in 2D,
the 3D algorithm estimates the normal vector with iteration algorithm.
Note that in this study, the reconstruction algorithm in this study follows
the algorithm of  \citet{lemoine_moment--fluid_2017} (in 2D)
and the algorithm of \citet{milcent_moment--fluid_2020} (in 3D).

Extended from VOF advection \citep{bna_review_2013},
the MOF advection can also be implemented either with operator split
\citep{jemison_coupled_2013,jemison_compressible_2014,asuri_mukundan_3d_2020}
or unsplit algorithm
\citep{dyadechko_reconstruction_2008,galera_2d_2011,breil_multi-material_2013},
with an additional advection step of the centroid.
\citet{dyadechko_moment--fluid_2005} show that
the material centroid in the incompressible flow move similarly to Lagrangian
particles and makes it a prefect choice for interface reconstruction input.
Unsplit methods have the advantage of only requiring one advection
and reconstruction step per time step,
the advection step is often algorithmically complex \citet{scardovelli_interface_2003}.
Split methods takes the advantage of the simple geometry of the Cartesian grid,
making the flux calculation simple and efficient.
Some simple and efficient formulation for the polyhedron intersection on
Cartesian has been proposed \citet{scardovelli_analytical_2000,milcent_moment--fluid_2020}.

There are some known issues on split advection algorithms,
especially on the mass conservation \citep{bna_review_2013}.
For a divergence-free flow field,
the discretized velocity field at intermediate time step as not
necessarily divergence-free.
The divergence correction term \citep{rider_reconstructing_1998} helps
to improve the local and global mass conservation.
The divergence correction term is expressed either in geometrical
for non-geometrical approach.
\citet{aulisa_geometrical_2003} show the geometrical properties of the
flux-correction term for Eulerian Implicit (EI) and Lagrangian Explicit (LE)
schemes and suggest using a hybrid EI and LE scheme helps improve the accuracy
and mass conservation.
The geometrical approach can be naturally extended to split MOF advection,
as the centroid can be determined directly from the corresponding geometrical
expression \citet{jemison_coupled_2013}.
The hybrid EI and LE scheme in 2D (EILE2D) \citep{aulisa_geometrical_2003} is
both geometrical preserving and area preserving,
and has been extended to MOF method
\citep{li_incompressible_2015,asuri_mukundan_3d_2020}.
The extension of EILE2D to 3D is not direct,
\citet{aulisa_interface_2007} composites the 3D velocity field
into 3 2D velocity fields.
however, this decomposition strategy
may cause stability issue in practical simulation.
Besides, the hybrid EI and LE scheme in 3D has not been extended to MOF advection.
Other approaches deals with the divergence-correction term with an algebaric
correction.
\citet{weymouth_conservative_2010} enforced the  divergence-correction term
with either 0 or 1,
depending on the initial value of the volume fraction.
This approach is ultimately mass conservation,
however, not geometrical preserving.
This approach has been extended to MOF advection,
with the centroids advection in a way that used in EI schemes
\citep{li_incompressible_2015,asuri_mukundan_3d_2020}.
Another non-geometrical correction introduces an Eulerian Algebraic (EA)
step between the EI and LE scheme (Namely EIEALE scheme) \citep{baraldi_mass-conserving_2014}.
Again, the intermediate Eulerian Algebraic does not guarantee geometrical
preserving
and has not yet been extended to MOF advection.

In this comparative study,
we investigates the accuracy and mass conservation properties of different split MOF
advection methods.
Besides the existing EI, LE, EILE2D \citep{aulisa_geometrical_2003,jemison_coupled_2013} and WY
schemes \citep{weymouth_conservative_2010,li_incompressible_2015,asuri_mukundan_3d_2020},
we also extend the EILE3D (Eulerian Implicit Lagrangian Explicit in 3D) \citep{aulisa_interface_2007},
, EILE3DS (EILE3D simplified) \citep{aulisa_interface_2007},
EIEALE (EI-EA-LE) scheme \citep{baraldi_mass-conserving_2014}
to from VOF-PLIC
advection to MOF advection.
None of the above methods guarantee both mass conservation and free stream
preserving for all scenarios.
We investigates the geometrical error, mass error and order of accuracy
with several 2D and 3D tests.
Besides,
we also proposed a simple and efficient algorithm that calculates
the volume fractions and centroids of the linear interface analytically.
This paper is organized as follows:
Section \ref{sec:govs} shows the goverining equations;
Section \ref{sec:MOFreconstruction} reviews the MOF reconstruction;
The implementation of the above-mentioned MOF advection schemes are introduced
in Section \ref{sec:MOFAdvection};
Section \ref{sec:geometrical} shows the geometrical algorithms for the linear
interface;
Several 2D and 3D numerical tests are investigated in Section \ref{sec:test};
Conclusions are drown in Section \ref{sec:conclusion}.

Note that the all results are based on the author's implementation,
The source code and test cases are available on https://github.com/zhoutengye/NNMOF.

\section{Governing Equations}
\label{sec:govs}
In MOF method, the characteristic function ( or color function) $c$ is used to indicates the fluid type
\begin{equation}
    \label{eq:colorfunction}
    c(\mathbf{x})=\left\{\begin{array}{ll}
        1 & \text { if } \mathbf{x} \in \text { fluid phase }  \\
        0 & \text { if } \mathbf{x} \in \text { other phases }
    \end{array}\right.
    .
\end{equation}
The volume fraction $C$ is a descretized version of the color function.
As an extension of the VOF-PLIC method,
the volume fraction is defined follows that in VOF method
\begin{equation}
    \label{eq:moffrac}
    C = \frac{\int_{\Omega}c(\mathbf{x})\rm{d}V}{\int_{\Omega}\rm{d}V},
\end{equation}
with the corresponding centroid defined from the first moment of the characteristic function
\begin{equation}
    \label{eq:mofcentroid}
    \mathbf{x_c} = \frac{\int_{\Omega} \mathbf{x_c}_c(\mathbf{x})\rm{d}V}{\int_{\Omega}c(\mathbf{x})\rm{d}V},
\end{equation}
where $\Omega$ is the domain of the descretized grid cell.
Giving a velocity field $\mathbf{u}$, a standard advection equation governs the evolution of
volume fraction $C$ and the centroid $\mathbf{x_c}$
\begin{equation}
    \begin{aligned}
        \label{eq:mofgov1}
        \frac{{\rm{d}}C}{dt} = 0 \rightarrow \frac{\partial C}{\partial t} + \mathbf{u} \cdot \nabla C = 0,
    \end{aligned}
\end{equation}
\begin{equation}
    \begin{aligned}
        \label{eq:mofgov2}
        \frac{{\rm{d}}\mathbf{x_c}}{dt} = \mathbf{u}.
    \end{aligned}
\end{equation}
For incompressible flow,
Eq. \eqref{eq:mofgov1} can be recast in conservation form
\begin{equation}
    \begin{aligned}
        \label{eq:mofgov1conservatinve}
        \frac{\partial C}{\partial t} + \nabla \cdot (C \mathbf{u}) = 0.
    \end{aligned}
\end{equation}
However,
in numerical descretization,
the velocity field may not always divergence-free,
especially for some multi-step schemes.
Eq. \eqref{eq:mofgov1conservatinve} is typically solved with a divergence correction term
\begin{equation}
    \begin{aligned}
        \label{eq:mofgov1correction}
        \frac{\partial C}{\partial t} + \nabla \cdot (C \mathbf{u}) = C \nabla \cdot \mathbf{u}.
    \end{aligned}
\end{equation}
\citet{rider_reconstructing_1998} show that even for incompressible flow,
the divergence correction term in Eq. \eqref{eq:mofgov1correction} help improving the local and global mass conservation.
In this study,
our algorithms are based on the divergence correction form.

The MOF method contains two key steps: (1) reconstruction algorithm (2) advection algorithm.
The reconstruction algorithm approximates the material in cell as a cut-off polyhedron (or cut-off polygon in 2D) from the computational grid cell,
the volume fraction and centroid of the cut-off in Eq. \eqref{eq:moffrac} is determined from
the cut-off polyhedron (polygon).
The advection algorithm solves Eq.  \eqref{eq:mofgov2} and Eq. \eqref{eq:mofgov1correction},
with the known information at step $n$ and advance to time step $n+1$.
Several important conditions should always satisfy for incompressible flow:

\begin{enumerate}
    \item Global conservation condition

          For incompressible flow, the global mass should always conserve.
          \begin{equation}
              \begin{aligned}
                  \label{eq:condition1}
                  \int_{\Omega_g} C^{n} \mathrm{d}V = \int_{\Omega_g} C^{n+1} \mathrm{d}V,
              \end{aligned}
          \end{equation}

    \item Local bound condition

          No overshooting or udershooting for volume fraction,
          which means no volume fraction smaller than the empty cell or greater
          than a full cell.
          In addition, the centroid should always be located within the grid
          cell $\Omega$
          \begin{equation}
              \begin{aligned}
                  \label{eq:condition2}
                   & 0 \le C \le 1           \\
                   & \mathbf{x_c} \in \Omega
              \end{aligned}.
          \end{equation}

    \item  Global geometrical consistency condition
          during an advection step,
          the sum of the pre-image or post-image
          should remain consistency,
          which means no hole or overlap should happen.

\end{enumerate}

The above 3 conditions are usually not be
satisfied simultaneously.
For split advection method,
even for a divergence-free velocity field,
the velocity field in each of the sub-step is not necessarily divergence free.
This may bring in additional bound violation for volume fraction and
the centroid.
the volume fraction is forced to be revalued to $\epsilon$ or $1-\epsilon$,
and the centroid is forced be revalued to $(\epsilon, \epsilon, \epsilon)$
or $(0.5,0.5,0.5)$.
Where $\epsilon$ is a very small number.
In this study, $\epsilon=1\times10^{-14}$.
If the correction happens,
the other two conditions, global conservation and global geometrical consistency
are no longer hold.
There are some way to enforce the global conservation or local bound conditions,
but the global geometrical consistency condition is rather stringent for
split advection method.

\section{MOF reconstruction}
\label{sec:MOFreconstruction}

The MOF method introduces the centroid as an additional constraint to determine
the reconstruction interface.
Unfortunately,
the known centroid $\mathbf{x_{c_{\rm{ref}}}}$ and volume fraction $C_{\rm{ref}}$
may not simultaneously satisfy with a linear cut-off.
The exact centroid matching is sacrificed on the altar of volume consercation.
A linear cut-off with the given volume fraction
which provides the best approximation to the reference centroid
is used as the optmimized cut-off for MOF method.

The MOF method reconstructs the interface in a 3D rectangular hextahedron cell
with a plane
\begin{equation}
    \label{eq:mofplane}
    \left\{\mathbf{x} \in \mathbb{R}^{3} \mid \mathbf{n} \cdot\left(\mathbf{x}-\mathbf{x}_0\right)+\alpha=0\right\},
\end{equation}
where $\mathbf{n}$ is the normal vector, $\mathbf{x_0}$ is  reference point of the cell,
either the center of the cell or the lower corner of the cell,
depending on the computational algorithm.
$b$ is the parameter that represents the distances from the reference point $\mathbf{x}_0$.
the volume of the cutting polyhedron by the reconstruction plane satisfies
\begin{equation}
    \label{eq:mofvolumeconstraint}
    \left|C_{\mathrm{ref}}-C_{A}(\boldsymbol{n}, b)\right|=0,
\end{equation}
where $C_\mathrm{ref}$ indicates the volume of the cutting polyhedron and
$C_A$ indicates the actual volume.
In addition to the the constraint on volume fraction,
the MOF reconstruction also minimizes error
between the reference centroid $\mathbf{x_{c_\mathrm{ref}}}$ and the actual centroid $\mathbf{x_{c_\mathrm{A}}}$
(See Fig. \ref{fig:sketchmof})
\begin{equation}
    \label{eq:mofcentroidconstraint}
    E_{\mathrm{MOF}}=\left\|\mathbf{x_{c_\mathrm{ref}}}-\mathbf{x_{c_\mathrm{A}}}(\mathbf{n}, b)\right\|_{2}.
\end{equation}
In 3D, Eq. \eqref{eq:mofcentroidconstraint} has 4 parameters.
As $b$ can be determined by the flood algorithm and the normal vector can be parameterized with
\begin{equation}
    \label{eq:norm2angle}
    \mathbf{n}=
    \left(\begin{array}{c}
            \sin (\phi) \cos (\theta) \\
            \sin (\phi) \sin (\theta) \\
            \cos (\phi)
        \end{array}\right),
\end{equation}
Eq. \eqref{eq:mofcentroidconstraint} is simplified as a function of $\phi$ and $\theta$.
\begin{equation}
    \label{eq:mofminimizeangle}
    E_{\mathrm{MOF}}\left(\phi^{*}, \theta^{*}\right)=\left\|\mathbf{f}\left(\phi^{*}, \theta^{*}\right)\right\|_{2}=\min _{(\phi, \theta):\text{Eq. (11) holds }}\|\mathbf{f}(\phi, \theta)\|_{2},
\end{equation}
where,
\begin{equation}
    \label{eq:parafunc}
    \mathbf{f}: \mathbb{R}^{2} \rightarrow \mathbb{R}^{3}, \quad \mathbf{f}(\phi, \theta)=\left(\mathbf{x}_{\text {ref }}-\mathbf{x}_{A}(\phi, \theta)\right).
\end{equation}
The minimization problem in Eq. \eqref{eq:mofminimizeangle} is a non-linear least square problem
for $\phi$ and $\theta$, which is solved numerically with optimization
algorithm.
Note that the expression for 2D problem is similar.

\ExecuteMetaData[figures.tex]{fig:sketchmof}

In the original MOF method \cite{dyadechko_moment--fluid_2005},
the non-linear optimization Eq. \eqref{eq:mofminimizeangle} is solved with
Broyden-Fletcher-Goldfarb-Shanno (BFGS) algorithm.
Several improvements has been done to accelerate the MOF reconstruction,
more specifically, on the geometrical algorithm to find the centroid during iteration
\citep{chen_improved_2016,milcent_moment--fluid_2020}.
\citet{lemoine_moment--fluid_2017} proposed an analytic algorithm
that calculates the optimized centroid directly,
without any numerical iteration.
Unfortunately,
this algorithm cannot be extended to 3D.

In this study, the analytic approach of \citet{lemoine_moment--fluid_2017} is used for 2D reconstruction,
and the numerical approach of \citet{milcent_moment--fluid_2020} is used for 3D reconstruction.

\section{MOF advection}
\label{sec:MOFAdvection}


Similar to VOF-PLIC method,
the advection of volume fraction and centroid in MOF method can be solved either with unsplit scheme  \cite{dyadechko_moment--fluid_2005}
and directional splitting scheme \cite{jemison_coupled_2013, asuri_mukundan_3d_2020}.
In either scheme,
the advection procedure is performed as mapping the fluid from from the departure region to the target region.
By finding the intersection between the fluid and the target region,
the volume fraction and centroids at new time step are determined.
It is important to guarantee the centroid associated to its phase packet during advection.
By using the fact that the material centroids in the incompressible flow move similarly to Lagrangian particles ( See proof in \cite{dyadechko_reconstruction_2008} Appendix A),
the consistency between Eq. \eqref{eq:mofgov2} and Eq. \eqref{eq:mofgov1correction} are guaranteed.
In this paper, we use directional splitting method for advection.
In rectangular hexahedron grid,
the sub-regions in the departure region and target regions are rectangular hexahedron.

Eq. \eqref{eq:mofgov1correction} is typically descretized with multi-step operator splitting algorithm,
a typicall descretization in $x-y-z$ coordinate contains three swppes:
\begin{equation}
    \begin{aligned}
        \label{eq:mof_adv_split3d}
        C_{i,j,k}^{(1)} & =C_{i,j,k}^{n} +\left(F_{i-1/2}^{u} -F_{i+1/2}^{u} \right)+C_{i, j, k}^{<1>}\left(u_{i+1 / 2}-u_{i-1 / 2}\right) ,  \\
        C_{i,j,k}^{(2)} & =C_{i,j,k}^{(1)} +\left(F_{j-1/2}^{v} -F_{j+1/2}^{v} \right)+C_{i, j, k}^{<2>}\left(v_{j+1 / 2}-v_{j-1 / 2}\right), \\
        C_{i,j,k}^{n+1} & =C_{i,j,k}^{(2)}+\left(F_{k-1/2}^{w}-F_{k+1/2}^{w}\right)+C_{i, j, k}^{<3>}\left(w_{k+1/2}-w_{k-1 / 2}\right).
    \end{aligned}
\end{equation}
Where $C^{(1)}_{i,j,k}$ and $C^{(2)}_{i,j,k}$ are the volume fraction at the intermediate step.
$F$ is the flux, will be determined either in an Lagrangian or Eulerian approach.
$C^{<>}$ means the volume function used for divergence correction,
the value is determined based on the advection scheme.
Note that the sequence of the sweep could change based on the implementation.
The volume fraction and flux in Eq. \eqref{eq:mof_adv_split3d} are determined from the intersection of the fluid and the
target region,
the centroid of the corresponding intersection can be determined using the scheme in the previous section.
In every directional split advection,
the new centroid can be calculated by weighted averaging
\begin{equation}
    \label{eq:mof_centroid}
    \mathbf{x_c}= \mathcal{M}\left({\frac{\sum V_{i} \mathbf{x}_i}{\sum V_{i}}}\right),
\end{equation}
where $i$ indicates the sub-region in the target region and $V_i$ indicates the volume fraction
$C$ or $F$ in Eq. \eqref{eq:mof_adv_split3d}.
$\mathcal{M}$ is the mapping form for the centroid,
which is based on the implementation of the advection scheme.


\subsection{Eulerian implicit scheme}
\label{sec:EI}

\ExecuteMetaData[figures.tex]{fig:EILEadvection}

The basic features for Eulerian implicit advection is plotted in Fig. \ref{fig:EILE_advection} (a).
For simplicity and clarity, we assume the $z$-conponent of the normal vector $n_z=0$,
so that the reconstructed polyhedron can be mapped to the $x-y$ plane as a 2D view.
We also assume the same velocity magnitude in all cell faces and the two boudnary velocity equals to 0.
The 5 mappings fron $n$ to $n+1$ in Fig. \ref{fig:EILE_advection} corresponds with
one-way expansion, linear translation, two-way compression, two-way expansion and one-way compression.

The key steps in Eulerian implicit advection scheme in $x$ direction are:

\begin{itemize}
    \begin{item}
          Determine the departure region and its sub-regions.

          By employing the upwinding $y-z$ plane (2D view in Fig. \ref{fig:EILE_advection} vertical line),
          every cell region $\Omega_{i}^{n}$ is now divided into smaller sub-regions.
          In directional splitting Eulerian implicit scheme, the mapping from the departure region to the target region of grid cell $i$ is defined as
          \begin{equation}
              \label{eq:ei_mapping}
              \Omega_{i,\rm{depart}}^n \rightarrow \Omega_{i,\rm{target}}^{n+1} =
              \left\{x_{i-1 / 2}-u_{i-1 / 2} \Delta t, x_{i+1 / 2}-u_{i+1 / 2} \Delta t\right\} \rightarrow\left\{x_{i-1 / 2}, x_{i+1 / 2}\right\}.
          \end{equation}
          The departure region $\Omega_{i}^{n}$ can be expressed as the union of two flux regions and one center region
          \begin{equation}
              \label{eq:ei_departure}
              \Omega_{i,\rm{depart}}^{n} = \Omega_{i,\rm{center}}^{n} \cup \Omega_{i-1/2,\rm{flux}}^{n} \cup \Omega_{i+1/2,\rm{flux}}^{n},
          \end{equation}
          where,
          \begin{equation}
              \label{eq:ei_fluxc}
              \Omega_{i,\rm{center}} = [ \max(x_{i-1/2},x_{i-1/2}-u_{i-1/2,}\Delta t), \min(x_{i+1/2}, x_{i+1/2}+u_{i+1/2,}\Delta t],
          \end{equation}
          \begin{equation}
              \begin{aligned}
                  \label{eq:ei_flux+}
                  \Omega_{i+1/2,\rm{flux}} & = \Omega_{i+1,\rm{target}}^n \cap \Omega_{i,\rm{depart}} \\
                                           & =
                  \left\{
                  \begin{aligned}
                      {\left[x_{i+1 / 2}, x_{i+1 / 2} + u_{i+1/2}\Delta t \right] }, &
                      \quad \quad \rm{if} \quad u_{i+1/2,} > 0                      \\
                      \emptyset,                                                     &
                      \quad \quad \rm{if} \quad u_{i+1/2,} < 0
                  \end{aligned}
                  \right .,
              \end{aligned}
          \end{equation}
          \begin{equation}
              \begin{aligned}
                  \label{eq:ei_flux-}
                  \Omega_{i-1/2,\rm{flux}} & = \Omega_{i-1,\rm{target}}^n \cap \Omega_{i,\rm{depart}} \\
                                           & =
                  \left\{
                  \begin{aligned}
                      {\left[x_{i-1 / 2}- u_{i-1/2}\Delta t, x_{i-1 / 2}  \right] }, &
                      \quad \quad \rm{if} \quad u_{i-1/2,} > 0                      \\
                      \emptyset,                                                     &
                      \quad \quad \rm{if} \quad u_{i-1/2,} < 0
                  \end{aligned}
                  \right . .                                                                            
              \end{aligned}
          \end{equation}

    \end{item}

    \begin{item}
          Find the volume fraction and the centroid of each sub-region using the reconstruction algorithm in the previous section.
    \end{item}

    \begin{item}

          Mapping the volume fraction and the centroid from the departure region to the target region

          We define the mapping ratio in Eulerian implicit scheme as
          \begin{equation}
              \label{eq:ei_compression}
              \beta = \frac{1}{1+u_{i-1/2,}\Delta t / \Delta x - u_{i+1/2,} \Delta t / \Delta x}.
          \end{equation}
          The value of $\beta$ indicates
          \begin{equation}
              \label{eq:ei_factor}
              \left\{\begin{array}{ll}
                  \beta<1, & \text { compression } \\
                  \beta=1, & \text { translation } \\
                  \beta>1, & \text { expansion }\end{array}\right. .
          \end{equation}

          The new volume fraction and centroids can be calculated with
          \begin{equation}
              \label{eq:ei_newphi}
              C_{i}^{*} = \beta\left ({V^{n}_{i,\rm{center}} + V_{i-1/2,\rm{flux}}^n + V_{i+1/2,\rm{flux}}} \right),
          \end{equation}
          \begin{equation}
              \label{eq:ei_newcen}
              \mathbf{x_c}_{i}^{*} = \beta\frac{
              \mathbf{x_c}^{n}_{i,\rm{center}} V^{n}_{i,\rm{center}}+
              \mathbf{x_c}_{i-1/2,\rm{flux}}^n V^{n}_{i-1/2,\rm{flux}}+
              \mathbf{x_c}^n_{i+1/2,\rm{flux}} V^{n}_{i+1/2,\rm{flux}}
              }{C^{*}_i} - \mathbf{x_0},
          \end{equation}
          where $\mathbf{x_0} = (x_i-1/2,y_j,z_k)$.

    \end{item}
\end{itemize}

Eq. \eqref{eq:ei_newphi} is equivalent to Eq. \eqref{eq:mof_adv_split3d}
when  the divergence correction term in $C^{<1>_{i,j,k}}=C^{n}_{i,j,k}$ in Eq. \eqref{eq:mof_adv_split3d}.
The EI mapping ensures all polyhedron in sub-regions aligned to the cell,
thus no overshooting or undershooting will happen,
the corresponding centroid is confined within the grid region,
thus the local bound condition holds during every split EI advection step.

\subsection{Lagrangian Explicit scheme}
\label{sec:LE}

The mapping strategy in Langantian Explicit scheme (LE) is slightly different from the EI scheme.
The departure region in LE scheme is $\Omega_i^n$ and the target region are determined by the
downwinding $y-z$ plane (2D view in Fig. \ref{fig:EILE_advection} vertical line).
In LE scheme,
the flux are calculated after mapping.
The key steps in LE advection are:

\begin{itemize}
    \begin{item}
          Determine the target region and  its sub-regions

          In directional splitting LE scheme, the mapping from the departure region to the target region of grid cell $i$ is defined as
          \begin{equation}
              \label{eq:le_mapping}
              \Omega_{i,\rm{depart}}^n \rightarrow \Omega_{i,\rm{target}}^{n+1} =
              \left\{x_{i-1 / 2}, x_{i+1 / 2}\right\}
              \rightarrow
              \left\{x_{i-1 / 2}+u_{i-1 / 2} \Delta t, x_{i+1 / 2}+u_{i+1 / 2} \Delta t\right\}.
          \end{equation}
          The target region $\Omega_{i}^{n}$ can be expressed as the union of two flux regions and one center region
          \begin{equation}
              \label{eq:le_target}
              \Omega_{i,\rm{target}}^{n+1} = \Omega_{i,\rm{center}}^{n+1} \cup \Omega_{i-1/2,\rm{flux}}^{n+1} \cup \Omega_{i+1/2,\rm{flux}}^{n+1},
          \end{equation}
          where,
          \begin{equation}
              \label{eq:le_targetc}
              \Omega_{i,\rm{center}}^{n+1} = [ \min(x_{i-1/2},x_{i-1/2}+u_{i-1/2,}\Delta t), \max(x_{i+1/2}, x_{i+1/2}+u_{i+1/2,}\Delta t],
          \end{equation}
          \begin{equation}
              \begin{aligned}
                  \label{eq:le_flux+}
                  \Omega_{i+1/2,\rm{flux}}^{n+1} & = \Omega_{i+1,\rm{depart}}^n \cap \Omega_{i,\rm{target}}^{n+1} \\
                                                 & =
                  \left\{
                  \begin{aligned}
                      {\left[x_{i+1 / 2}, x_{i+1 / 2} + u_{i+1/2}\Delta t \right] }, &
                      \quad \quad \rm{if} \quad u_{i+1/2,} > 0                      \\
                      \emptyset,                                                     &
                      \quad \quad \rm{if} \quad u_{i+1/2,} > 0
                  \end{aligned}
                  \right . ,
              \end{aligned}
          \end{equation}
          \begin{equation}
              \begin{aligned}
                  \label{eq:le_flux-}
                  \Omega_{i-1/2,\rm{flux}}^{n+1} & = \Omega_{i-1,\rm{depart}}^n \cap \Omega_{i,\rm{target}}^{n+1} \\
                                                 & =
                  \left\{
                  \begin{aligned}
                      {\left[x_{i-1 / 2}+ u_{i-1/2}\Delta t, x_{i-1 / 2}  \right] }, &
                      \quad \quad \rm{if} \quad u_{i-1/2,} > 0                      \\
                      \emptyset,                                                     &
                      \quad \quad \rm{if} \quad u_{i-1/2,} < 0
                  \end{aligned}
                  \right .  .                                                                                      \\
              \end{aligned}
          \end{equation}

    \end{item}

    \begin{item}
          Mapping the volume fraction and the centroid from the departure region to the target region.

          In LE scheme,
          The flux is not determined from the departure region,
          but is estimated from the intersection of remapped polyhedron and the target region and the grid cells.
          The mapping factor in LE advection is defined as
          \begin{equation}
              \label{eq:le_compression}
              \gamma = 1-u_{i-1/2,}\Delta t / \Delta x + u_{i+1/2,} \Delta t / \Delta x.
          \end{equation}
          Similar to the expansion factor for the EI advection in Eq. \eqref{eq:ei_compression},
          the value of $\gamma$ indicates
          \begin{equation}
              \label{eq:le_factor}
              \left\{\begin{array}{ll}
                  \gamma<1, & \text { compression } \\
                  \gamma=1, & \text { translation } \\
                  \gamma>1, & \text { expansion }\end{array}\right. .
          \end{equation}
          After mapping, the normal vector of the reconstruction plane changes,
          the new normal vector can be determined by altering the value $(\alpha, m_x)$ with
          \begin{equation}
              \label{eq:le_shift}
              \begin{aligned}
                  m_x^\prime & = \gamma m_x,                         \\
                  \alpha'    & = \alpha + m_x u \Delta t / \Delta x.
              \end{aligned}
          \end{equation}
    \end{item}

    \begin{item}

          Calculate the flux and center region and update volume fraction and centroids

          The new volume fraction and centroid in the sub-regions can be determined by using the algorithm in previous section.
          The volume fraction and centroids in $\Omega^{n+1}_{i}$ can be determined by

          The new volume fraction and volume fraction are mapped from the
          \begin{equation}
              \label{eq:le_newphi}
              C^{*} = \left ({V^{n}_{i,\rm{center}} + V_{i-1/2,\rm{flux}}^n + V_{i+1/2,\rm{flux}}} \right),
          \end{equation}
          \begin{equation}
              \label{eq:le_newcen}
              \mathbf{x_c}^{*} = \frac{
              \mathbf{x_c}^{n}_{i,\rm{center}} V^{n}_{i,\rm{center}}+
              \mathbf{x_c}_{i-1/2,\rm{flux}}^n V^{n}_{i-1/2,\rm{flux}}+
              \mathbf{x_c}_{i+1/2,\rm{flux}} V^{n}_{i+1/2,\rm{flux}}
              }{C^{*}_i}.
          \end{equation}

    \end{item}
\end{itemize}

The implementation of LE scheme is simple.
Although the the LE advection can be expressed as a mapping procedure (See Fig. \ref{fig:EILE_advection}),
the mapping procedure does not have to be solved explicitly.
Simply by altering the slope the $\alpha$,
the sub-regions in $\Omega^{n+1}$ are determined.
When the divergence correction term in $C^{<1>_{i,j,k}}=C^{n}_{i,j,k}$ in Eq.
\eqref{eq:mof_adv_split3d},
the expression is equivalent to Eq. \ref{eq:le_newphi}.
LE advection also automatically takes the divergence term into account.
Same as EI advection,
a single split LE advection satisfies the local bound condition for both volume
fraction and centroid.

\subsection{EI-LE strategy in 2D}
\label{sec:EILE2D}

\ExecuteMetaData[figures.tex]{fig:eile2d}

\citet{scardovelli_interface_2003} show that both EI and LE operator splitting advection schemes
satisfies the local bound condition,
neither of them guarantees the global mass conservation if only one of them is used.
The main reason that causes the unconserved mass is the overlap
or "hole" of between the departure and arrival regions between two time steps,
which consists of multiple split advections.
The global geometrical consistency condition will not hold as well,
which brings in additional error to the advection of centroid.
A hybrid EI-LE scheme \citep{scardovelli_interface_2003,aulisa_geometrical_2003}
(in this paper, referred to as EILE2D)
is proposed which ultimately conserves volume.
Assume the operator splitting advection starts with an EI advection in $x$-direction and
subsequently applies an LE advection $y$-direction (See Fig. \ref{fig:eile2d} (a)),
the two mapping steps are
\begin{equation}
    \begin{aligned}
        \label{eq:eilemapping}
        \Pi_x = & \Omega_{\rm{depart}} \rightarrow \Omega_{i,j} =
        \left\{\begin{array}{l}
            x=\beta_x \left(x^\prime+u_{1}\right) \\
            y=y^\prime
        \end{array}\right. ,                    \\
        \Pi_y = & \Omega_{i,j}^n \rightarrow \Omega_{\rm{target}} =
        \left\{\begin{array}{l}
            x^{\prime\prime}= x \\
            y^{\prime\prime}=\gamma_y y+v_{1}
        \end{array}\right.  ,                  \\
    \end{aligned}
\end{equation}
where,
\begin{equation}
    \begin{aligned}
        \label{eq:betagamma}
        \beta_x = \frac{1}{1+u_1 - u_2}, \\
        \gamma_y = 1-v_1 + v_2 .          \\
    \end{aligned}
\end{equation}
The mapping from the departure region and target region $\Pi_{xy}$ can be expressed as
\begin{equation}
    \begin{aligned}
        \label{eq:eilemapping2}
        \Pi_{xy} = \Pi_{x}\Pi_{y}=
        \left\{\begin{array}{l}
            x^{\prime\prime}=\beta_x \left(x^\prime+u_{1}\right) \\
            y^{\prime\prime}=\gamma_y y^\prime + v_1 
        \end{array}\right. \\
    \end{aligned} .
\end{equation}
With the divergence-free constraint
\begin{equation}
    \label{eq:divcond}
    u_2 - u_1 = v_1 - v_2.
\end{equation}
The Jacobian of the transformation $J=\beta_x \gamma_y=1$.
This transformation is an volume-preserving mapping.
When the advection starts with $y$-direction, the
EI-LE scheme also guarantees volume-preserving (See Fig. \ref{fig:eile2d} (b)),
thus condition Eq. \eqref{eq:condition2} satisfies.
Since both EI and LE advection are consistent with condition Eq. \eqref{eq:condition1},
the EI-LE scheme satisfied both of the conditions required by VOF.
Note that the mapping factors in Eq. \eqref{eq:eilemapping2} are consistent
with the mapping factors in Eq. \eqref{eq:ei_compression} and Eq. \eqref{eq:le_compression},
the volume-preserving also guarantees the correct mapping of the centroid during
MOF advection.
If the sweep is done with a reversed order,
which starts with an LE advection and then an EI advection (referred to as
LEEI2D),
none of global conservation or global geometrical consistency condition will hold.

\subsection{EI-LE strategy in 3D}
\label{sec:EILE3D}
Unfortunately,
the EILE strategy cannot extend to 3D directly.
\citet{aulisa_geometrical_2003} show that
the ultimate geometrical preserving in 3D can be achieved by
decompositing the velocity field to three 2D
divergence-free velocity fields
$\mathbf{u_1} = (\bar u_1,\bar v_1,0)$,
$\mathbf{u_2} = (\bar u_2,0,\bar w_2)$,
$\mathbf{u_3} = (0,\bar v_3,\bar w_3)$
which satisfied
\begin{equation}
    \begin{aligned}
        \label{eq:eile3ddecomposition}
         & \mathbf{u} = \mathbf{u_1} + \mathbf{u_2} + \mathbf{u_3}, \\
         & \nabla \cdot \mathbf{u_1} = 0,                           \\
         & \nabla \cdot \mathbf{u_2} = 0,                           \\
         & \nabla \cdot \mathbf{u_3} = 0.
    \end{aligned}
\end{equation}
One simple approach
to composite the 3D velocity field is to
use the fact that the only 5 of the 6 unknowns
$(\bar u_1, \bar u_2, \bar v_1, \bar v_2, \bar w_1, \bar w_2)$
are independent.
By set $\bar u_1 = u/2$ for the whole domain,
and the lower boundary value of
$\bar v_{1_{i,j-1/2,k}} = v_{{i,j-1/2,k}}/2$,
solve for $\bar v_{1_{i,j+1/2,k}}$
using the relationship $\nabla \cdot \mathbf{u_1}=0$.
The other components can be solved recursively by using
other divergence-free constraints.
Although the decomposition starts with a halved velocity component,
the velocity decomposition may not guarantee the magnitude of the velocity
bounded to the maximum magnitude of the original velocity field.
An example of unbounded decomposition is shown in Section \ref{sec:sv3d}.

A simplified EI-LE strategy in 3D (EILE3DS) is to alter the EI and LE sweeps
in two time steps.
For odd step, the sweep is X(EI)-Y(LE)-Z(EI);
while for even step, sweep is X(LE)-Y(EI)-Z(LE) \citep{li_incompressible_2015}.
This approach does not guarantee volume preserving,
however, improved the mass conservation of pure EI or LE advection.
Note that the strategy \citet{aulisa_interface_2007} also named EILE3DS,
however, implemented with 6 split advection for a full time step,
which is different implementation from this study.

\subsection{EI-EA-LE strategy in 3D}
\label{sec:EIEALE}
In this study, we proposed a simplified EI-LE splitting scheme for
both centroid and volume fraction.
\citet{baraldi_mass-conserving_2014} introduced an EA (Eulerian algebaric) step
between the EI and LE advection.
The mapping at EA step is determined by enforcing the Jacobian $J=1$ for the mapping
$\Pi_{xyz} = \Omega_{\rm{depart}} \rightarrow \Omega_{\rm{target}}$,
which does not have a geometrical interpretation.
With the approximate intermediate step,
the mass conservation in significant improved compared with
EI or LE advection.
In the numerical test of \citet{baraldi_mass-conserving_2014},
the mass conservation reaches machine precision with their
wisp-suppression algorithm.
In this study,
we use a similar approach for both volume fraction and centroid advection,
without wisp-suppression algorithm.

The advection consists of three steps.
Consider the mapping $\Pi_{xyz}$ corresponds with the sweep sequence of $x-y-z$,
the mapping during the advection from $n$ to $n+1$ is
\begin{equation}
    \begin{aligned}
        \label{eq:eiealemapping}
        \Pi_{xyz} = & \Omega_{\rm{depart}} \rightarrow \Omega_{\rm{target}} =
        \left\{\begin{array}{l}
            x=\beta_x \left(x^\prime+u_{1}\right)          \\
            y=\xi_{y1} \left(\xi_{y2}y^\prime + v_1\right) \\
            z=\gamma_z z^\prime +w_1                       \\
        \end{array}\right.                              \\
    \end{aligned}
\end{equation}
The expression of $\beta_x$ and $\gamma_z$ are consistent with
Eq. \eqref{eq:betagamma} in EILE2D scheme.
The intermediate step is an algebraic step,
to ensure that the Jacobian $\beta_x \xi_{y1} \xi_{y2} \gamma_z$ equal to 1,
$\xi_{y1}$ and $\xi_{y2}$ are expressed as
\begin{equation}
    \begin{aligned}
        \label{eq:xi1xi2}
        \xi_{y1} = \frac{1}{\gamma_z}, \quad \xi_{y2} = \frac{1}{\beta_x}.
    \end{aligned}
\end{equation}
For split advection procedure, the mapping Eq. \eqref{eq:eiealemapping} is
implemented with
\begin{equation}
    \label{eq:eiealemapping2}
    \begin{aligned}
         & \Pi_x=\Omega_{\rm{depart_1}} \rightarrow \Omega_{\rm{target_1}} = x=\beta_x \left(x^\prime+u_{1}\right),          \\
         & \Pi_y=\Omega_{\rm{depart_2}} \rightarrow \Omega_{\rm{target_2}} = y=\xi_{y1} \left(\xi_{y2}y^\prime + v_1\right), \\
         & \Pi_z=\Omega_{\rm{depart_3}} \rightarrow \Omega_{\rm{target_3}} = z=\gamma_z z^\prime +w_1 .                      \\
    \end{aligned}
\end{equation}
Because of the non-geometrical property of the intermediate EA step,
$\Omega_{\rm{target_1}} \ne \Omega_{\rm{depart_2}}$ and
$\Omega_{\rm{target_2}} \ne \Omega_{\rm{depart_3}}$,
the mapping $\Pi_{xyz} \ne \Pi_x \Pi_y \Pi_z $.
Altough the Jacobian of $\Pi_{xyz}$ equals to 1,
the Jacobian of $\Pi_{x}\Pi_{y}\Pi_{z}$ would not always equal to 1.

The expression of the volume fraction and centroid of the $x$ sweep and $z$
sweep
corresponds with the expression in EI and LE advection, respectively.
The volume fraction in EA step are calculated with
\begin{equation}
    \begin{aligned}
        \label{eq:eavol}
        C_{i j, k}^{(2)}=\frac{C_{i j, k}^{(1)}\left(1-\delta_{x} \bar{u}\right)+F_{j-1}^{v}-F_{j}^{v}}{1+\delta_{z} \bar{w}}.
    \end{aligned}
\end{equation}
For the centroid,
since there is no corresponding geometrical mapping for the EA step,
the centroid is calculated with the EI method as in Eq. \eqref{eq:ei_newcen}.

We have not found any proof that the EA step satisfies the local bound
condition due to the non-geometrical property.
Numerical test show that EIEILA scheme would not guarantee the ultimate mass conservation.

\subsection{WY approach}
\label{sec:WY}
The strategy by \citet{weymouth_conservative_2010} can be categorized as an EI advection,
as the flux calculation follows the EI approach.
The flux correction term $C^{<>}$ in Eq. \eqref{eq:mof_adv_split3d} is applied based on the value of the volume fraction at step $n$
\begin{equation}
    \label{eq:wycorrection}
    C^{<1>}_{i,j,k}(u_{i+1/2}-u_{i-1/2}) = \tilde C_{i,j,k} = \left\{
    \begin{array}{l}
        0,  \quad C^{n}_{i,j,k}>=\frac{1}{2} \\
        1,  \quad C^{n}_{i,j,k}<\frac{1}{2}  \\
    \end{array}
    \right. .
\end{equation}
The correction does not rely on the velocity field,
but is only based on the volume fraction at time step $n$.
\citet{weymouth_conservative_2010} prove that with the threshold $1/2$,
there is no over shooting or undershooting for the volume fraction,
and ultimately mass preserving.
However,
the flux correct term is an algebraic correction,
there is no corresponding geometrical expression.
There is no corresponding form for the correction in Eq. \eqref{eq:wycorrection}
for the centroid.
The EI centroid advection is applied along with the WY volume correction in
\citet{li_incompressible_2015}.
In this study,
we use the same strategy to the centroid.
The EI advection of the centroid guarantees the centroid in the grid cell.
With the volume fraction bounded between 0 and 1,
the local bound condition of MOF advection holds.
Due to the non-geometrical correction,
the global geometrical consistency condition does not hold for WY approach.

\subsection{Summarize}
\label{sec:mofadvcompare}

Among all fluxed-splitting MOF advection methods mentioned
in this section,
the key idea is to manipulate the divergence correction term,
either by changing the mapping strategy or apply modification to correction term directly.
The 2D and 3D methods are summarized in
Table \ref{tab:comparemethods2d} and
Table \ref{tab:comparemethods3d}.

\ExecuteMetaData[tables.tex]{tab:comparemethods2d}

\ExecuteMetaData[tables.tex]{tab:comparemethods3d}

For 2D methods,
all methods satisfies the local bounded condition,
which means the mass conservation issue does not
come from the correction from overshooting or undershooting,
but comes from the geometrical inconsistency.
EI nor LE would not conserve mass because its
non-volume preserving mapping.
For EILE2D method, with an LE advection coming after an EI advection,
satisfies all three conditions.
While the LEEI2D method,
although alters between LE and EI schemes at two sub-steps,
would not guarantee geometrical preserving or geometrical consistency.
The WY scheme ultimately conserve mass with a more stringent
CFL limit than other schemes,
however, due to a lack of the divergence correction,
failed on the geometrical consistency condition.
All the 2D strategies takes 2 sweeps in one time step.

For 3D methods,
pure EI or LE scheme involves mass conservation issue as
same properties to its corresponding 2D schemes.
For EILE3D scheme,
if the velocity decomposition guarantees divergence-free,
all three conditions hold.
The CFL limit seems to be 1,
as the starting values of the velocity field is
a halved velocity component,
however, the decomposited velocity may not be bounded to a half of the velocity
field
the CFL limit may not be one, or could be even much smaller than 1/2.
Besides, it takes twice of the advection steps of other advection methods.
The EILE3DS alters between EI and LE schemes,
is shown to preserve mass better than pure EI or LE scheme
\citep{jemison_compressible_2014}.
For EIEALE scheme,
the mapping coefficient of the EA step is determined from the Jacobian equal to
1,
although would not guarantee ultimate mass conservation,
helps to improve the mass conservation \citep{baraldi_mass-conserving_2014}.
While WY scheme in 3D has a even more stringent CFL limit than the 2D WY scheme,
keeping mass conservation within its CFL limit of $1/6$.
Note that all methods,
the local bounded condition holds for all methods,
excepts for EIEALE method.
Even though the EIEALE significant improve the mass conservation from EI or LE
method (See Section \ref{sec:test}),
we have not find any proof to show that the EA step guarantee the local bound
condition for volume fraction and centroid.

\section{Geometrical algorithm for MOF}
\label{sec:geometrical}

In MOF method, two problems are typically solved by geometrical algorithm: \\
(1) Intersection algorithm: Find the centroid $\mathbf{x_c}$ and volume fraction $C$ of the
cut-off polyhedra by a known plane (or in 2D the cutting polyhedra by a line).  \\
(2) Flood algorithm: Find the centroid $\mathbf{x_c}$ for the cutting plane (or line in 2D) with known normal vector $\mathbf{n}$ and volume fraction $C_{\rm{ref}}$.

During the MOF reconstruction,
the flood algorithm is typically used to determine the equation of the plane (line)
and the corresponding centroid of the cutting polyhedron (polygon)
with the normal vector in the current step.
The MOF advection,
in which is equation of the plane (line) is known,
the intersection algorithm is used to determine the
volume fraction and centroid in departure and target regions.
Note that most of the flood algorithms
contains the intersection algorithm as a part of it.

Many polyhedron (polygon) algorithms/tools have been developed
for the intersection and flood algorithms for VOF-PLIC method or MOF method
\citep{
    nijenhuis_combinatorial_2014,
    ahn_geometric_2008,
    rider_reconstructing_1998,
    lopez_analytical_2008,
    diot_interface_2014,
    diot_interface_2016}.
Anthough those algorithms works for general grids,
they have to deal with the intersection and the volume/centroid calculation
using a complex node/edge/face algorithm.
For some simple shapes,
analytic formulation can be found to accelerate the volume and centroid calculation
\citep{scardovelli_analytical_2000,yang_analytic_2006}.

In this study,
we implement very simple and efficient intersection and flood algorithms.
With the fact that all cell and flux regions during reconstruction and
advection are rectangular hexahedron (or rectangle),
the intersection becomes much simpler.
Unlike the computational geometrical approach,
the operators for the vertex, edge or face does not have to be calculated
explicitly. Here we only show the 3D algorithms,
the implementation of the 2D algorithm is similar,
but simpler.

\ExecuteMetaData[figures.tex]{fig:cuts}

The volume fraction in the intersection algorithm follows
\citet{scardovelli_analytical_2000} which
derived is derived with the combination of right tetrahedrons.
By implementation symmetric, permutation and shift origin techniques,
the cutting polyhedron can be simplified to 5 configurations (See Fig. \ref{fig:cuts}).
We derive the expression for the centroids as an extension of the approach of \citet{scardovelli_analytical_2000}.
The volume fraction
and the corresonding centroid are expressed as
\begin{equation}
    \label{eq:cuts}
    V=\frac{1}{6 m_{1} m_{2} m_{3}}\left[\alpha^{3}-\sum_{i=1}^{3} F_{3}\left(\alpha-m_{i} \Delta x_{i}\right)+\sum_{i=1}^{3} F_{3}\left(\alpha-\alpha_{\max }+m_{i} \Delta x_{i}\right)\right]
\end{equation}
\begin{equation}
    \begin{aligned}
        \label{eq:centroids}
        c_j=
        \frac{\alpha^4}{24 m_{j}V}
         & \left[
        - \sum_{i=1}^{3} (\alpha+(1+2\beta) m_i\Delta x_i) F_{3} (\alpha-m_{i} \Delta x_{i}) \right.                                           \\
         & \left . +  \sum_{i=1}^{3} (\alpha+(1-2\beta) (\alpha_{\max } - m_{i} \Delta x_{i})) F_{3}(\alpha-\alpha_{\max }+m_{i} \Delta x_{i})
        \right]
    \end{aligned}
\end{equation}
Where
\begin{equation}
    \label{eq:al}
    \alpha_{\max }=\sum_{i=1}^{3} m_{i} \Delta x_{i},
\end{equation}
\begin{equation}
    \label{eq:bt}
    \beta=\left\{\begin{array}{ll}
        -1 & \text { for } i \ne j \\
        1  & \text { for } i = j
    \end{array}\right.,
\end{equation}
\begin{equation}
    \label{eq:f3}
    F_{n}(y)=\left\{\begin{array}{ll}
        y^{n} & \text { for } y>0  \\
        0     & \text { for } y<=0
    \end{array}\right.,
\end{equation}

\ExecuteMetaData[figures.tex]{fig:flood}

Based on the values of $\mathbf{n}$ and $\alpha$,
the configuration is determined from the 5 candidate
cutting polyhedra and the volume fraction and centroid
are calculated using a very simple and efficiency equations.
The equation of the volume fraction and centroid for the 5
configurations are derived from Eq. \eqref{eq:cuts} and Eq. \eqref{eq:centroids}.
The detailed expressions are listed in Appendix A.

In flood algorithm,
the parameter $\alpha$ can be determined directly from the inverse function of Eq.
\eqref{eq:cuts},
the corresponding expressions for the 5 candidate configurations can be found in
\citet{scardovelli_analytical_2000}.
Although there is no explicit inverse function for Eq. \eqref{eq:centroids},
as the volume fraction is determined from the inverse form of Eq.
\eqref{eq:cuts},
the centroid can now be determined from the volume fraction and the expression
of the plane.

We applied both our algorithm and the geometrical algorithm of \citet{jemison_coupled_2013}
to the split advection methods in Section \ref{sec:MOFreconstruction},
on a 3D advection test (Corresponds with the Zalesak's test in Section \ref{sec:zalesak3d}).
The comparison of the CPU times is shown in Fig. \ref{tab:cputime}.
Overall, using the same advection and reconstruction,
our geometrical algorithm is about 50 times faster than the geometrical algorithm of \citet{jemison_coupled_2013}.
With the geometrical algorithm of \citet{jemison_coupled_2013},
the advection takes about 78\% of the computational time during advection.
While using our geometrical algorithm,
the ratio of CPU time the on advection term reduces from 78\% to 8\%.

\ExecuteMetaData[tables.tex]{tab:cputime}

\section{Numerical Tests}
\label{sec:test}

In this section,
six typical tests are taken to evaluate the performance
of different split MOF advection methods.
Three of them are in 2D,
and other three of the are in 3D.
The candidate methods for 2D and 3D are listed in Table
\ref{tab:comparemethods2d} and \ref{tab:comparemethods3d}.
For all the numerical tests,
the reconstruction algorithm, centroid algorithm and flood algorithm are identical.
For MOF reconstruction algorithm,
the maximum iteration is 10,
and the tolerance for iteration is $10^{-8}$.

Two different error measurement criteria \citep{zhang_new_2008} are used in this study:

(1) Geometrical error
\begin{equation}
    \label{eq:geometrical_error}
    E_{g}=\frac{\sum\left|f(T)-f_{exact}\right|}{N^{\mathrm{dim}}},
\end{equation}

(2) Mass error
\begin{equation}
    \label{eq:mass_error}
    E_{m}=\frac{\sum\left|f(T)-f_{\mathrm{exact}}\right|}{N^{\mathrm{dim}}},
\end{equation}
where $f(T)$ is the numerical results at $t=T$,
$f_{\rm{exact}}$ is the exact results,
$N$ is the grid resolution in $x$ direction,
$\rm{dim}=2$ for 2D problem and $\rm{dim}=3$ for 3D problem.
The order of the accuracy \cite{zhang_new_2008,aulisa_interface_2007} is defined as
\begin{equation}
    \label{eq:accuracy}
    \mathcal{O}_{h}=\log _{2}\left(\frac{E^{\rm{dim}}\left(\frac{1}{2 h}\right)}{E_{g}\left(\frac{1}{h}\right)}\right)
\end{equation}

Note that,
the CFL numer in the tests can be larger than the CFL limit listed
in Table \ref{tab:comparemethods2d} and Table \ref{tab:comparemethods3d}.

\subsection{2D tests}
\label{sec:2dtest}
\subsubsection{Zalesak's disk rotation test}
\label{sec:zalesak2d}

\ExecuteMetaData[figures.tex]{fig:2dzalesakinit}

The Zalesak's disk rotation test is firstly introduced by
\citet{zalesak_fully_1979}
and used by many other studies
\citep{rudman_volume-tracking_1997,aulisa_geometrical_2003,aulisa_interface_2007,zhang_new_2008}.
In Zalesak's disk rotation,
the rotation velocity field is defined with the a rotation velocity field.
The setup and velocity field of the Zalesak's rotation test is shown in Fig. \ref{fig:2dzalesakinit}.
After a full revolution of $2\pi$ rotation
the notched disk returns to its initial location.

\ExecuteMetaData[figures.tex]{fig:2dzalesakcfls}

\ExecuteMetaData[figures.tex]{fig:2dzalesakresos}

The geometrical errors are shown in Fig. \ref{fig:2dzalesakcfls} and Fig. \ref{fig:2dzalesakresos}.
For a fixed CFL number, the geometrical error decreases with the increase of the grid resolution.
While for a fixed grid resolution,
the geometrical error get larger when the CFL number becomes smaller.
In this test,
all methods do not show significant difference from each other when the CFL number
and grid resolution $N$ are fixed.
For mass conservation,
all methods conserve to machines precision with all combinations of
CFL number and grid resolution $N$.

\subsubsection{Rider-Kothe Single vortex test}
\label{sec:sv2d}

First introduced by \citet{rider_reconstructing_1998},
the single vortex test
represents the evolution of a circular disk under a
reversed shear velocity field.
The initial set up and test parameters are shown in Fig. \ref{fig:2dsvinit}.
With the domain-centered vortex field,
the material interface of the circular disk stretches into
a filament towards the center of the votrex.
The material interface reaches its maximum deformation
at $t=T/2$,
and recovers back to the initial circular disk at $t=T$.
The material interface may get teared or smoothed out
there is no sufficient numerial resolution.
Note that even though the stability criteria for most of the MOF advection methods is $\rm{CFL} \le 1/2$,
the $\rm{CFL}$ limit for this test could be greater then 0.5.
We test $\rm{CFL} = 0.8$ and $\rm{CFL} = 1.0$ as well.

\ExecuteMetaData[figures.tex]{fig:2dsvinit}

\ExecuteMetaData[figures.tex]{fig:2dsvcontour}

When the grid resolution is $64 \times 64$,
the material interfaces with different CFL number at
$t=T/2$ and $t=T$ are plotted in Fig. \ref{fig:2dsvcontour}.
When $t=T/2$,
all advection methods maintains the overall agreements with
the exact resolution,
with sharp corner on the tail and head smoothed out.
There is no significant difference among all advection methods.
When $t=T$,
the material interfaces show significant difference among each other,
especially when the CFL number is large.
Compared with the exact circular disk,
the EI and LE methods show most significant deviation.
The position of the EI material interface is left of the exact result,
while the position of the LE material interface right of the exact result.
Other methods,
although lost some accuracy at the top of the circular disk,
remains close to the exact material interface.

\ExecuteMetaData[figures.tex]{fig:2dsvresos}

The Geometrical error and mass conservation error with
grid resolution of $64 \times 64$ are plotted
in Fig. \ref{fig:2dsvresos}.
The optimal CFL number that results in minimum geometrical error for each method occurs at the range of $[0.05-0.2]$.
When CFL number is small ($<0.1$),
the geometrical error among all methods is not significant.
After $CFL>0.15$, with the increase of CFL number,
the geometrical error increases for most of the methods.
The EILE2D method and LEED2D method has overall smaller
geometrical error compared with EI, LE and WY methods.
For mass conservation,
EILE2D method conserves mass to machine precision as expected.
EI and LE has most significant mass error,
especially when CFL number is large.
When the CFL number is getting smaller,
the mass error becomes smaller.
When CFL number is 0.05,
the mass error is approximately $10^{-3}$.
The LEEI2D method which has a reversed order of EI and LE combination to EILE2D methods,
although does not satisfied the global mass conservation condition,
improves the mass error compared with EI and LE methods.
It is interesting to find that although the CFL number need to be smaller than $1/4$ in order to guarantee an ultimate mass preserving,
the WY method preserves mass even when CFL number is 1.0
in this study.

\ExecuteMetaData[figures.tex]{fig:2dsvcfls}

The geometrical error and mass error with the change of grid resolution are plotted in Fig. \ref{fig:2dsvcfls}.
The geometrical error decreases when the grid resolution
is getting larger.
When CFL number is small (CFL=0.05 as shown in Fig. \ref{fig:2dsvcfls} (a)),
the geometrical error are similar to each other for all methods.
When CfL number is large (CFL=1.0 as shown in Fig. \ref{fig:2dsvcfls} (b)),
EILE2D and LEEI2D shows a smaller geometrical error and higher order of accuracy.
The EILE2D and WY methods conserves the mass to machine
precision.
EILE2D method improves the mass conservation compared with EI or LE method.

\subsubsection{Reverse Vortex test}
\label{sec:rv2d}

We introduce a even more stringent test compared
with Rider-Kothe single vortex test \citep{rider_reconstructing_1998}.
The reverse vortex test is firstly introduced by
\citet{leveque_high-resolution_1996}.
The initial set up and test parameters are shown in Fig. \ref{fig:2drvinit}.
The number of vortex is set to 4 in this study,
which means a matrix of $4 \time 4$ vortex in the computational domain.
The domain-centered circular suffers a even more severe
deformation then that in Rider-Kothe single vortex test.
As shown in Fig. \ref{fig:2drvinit},
when $t=T/2$, although the material interface remains simply connected,
some of the filament of the interface are very thin.
When the period $T=2$,
the width of the thin filament can be as small as $2\times 10^{-4}$.
When the grid resolution is not enough,
the material interface will tear.
Similar to the Rider-Kothe single vortex test,
this deformation test contains the computational cases with
$\rm{CFL} = 0.8$ and $\rm{CFL} = 1.0$ as well.

\ExecuteMetaData[figures.tex]{fig:2drvinit}

\ExecuteMetaData[figures.tex]{fig:2drvcontour}

When the grid resolution is $64 \times 64$,
the material interfaces with different CFL number at
$t=T/2$ and $t=T$ are plotted in Fig. \ref{fig:2dsvcontour}.
When $t=T/2$,
teared material material interfaces are observed
for all methods.
However,
all methods remains good agreement with the main body of the material domain around center of the computational domain.
At the time $t=T$,
when CFL number is small (0.05),
all methods recover to the initial shape with reasonable
solution and do not show much difference between each other.
When the CFL number is larger,
EI and LE methods show a larger deviation with the
exact solution,
while the other three methods shows better agreement.

\ExecuteMetaData[figures.tex]{fig:2drvresos}

The geometrical error and mass conservation error with
grid resolution of $64 \times 64$ are plotted
in Fig. \ref{fig:2drvresos}.
The overall tendency of the the geometrical error with
the change of CFL number is to decrease when the
CFL number is getting smaller.
The EI methods has the overall maximum error
among all methods,
LE and WY methods comes after with similar magnitude of
of gemertric error.
The EILE2D method has similar performance to LE and
WY method when $CFL \le 0.25$,
and performs better when $CFL >0.25$.
The LEEI2D method has the overall minimum error
among all methods.
The EI method has the largest mass error,
and LE method comes after.
The LEEI2D method improves the mass conservation
of the EI and LE methods.
The EILE2D method preseves mass to machine precision
for all CFL numbers in the test cases.
The WY methods loses its ultimate mass preseving at $\rm{CFL}=0.8$ and $\rm{CFL}=1.0$.
Even for the caese when WY methods loses its ultimate mass conservation,
the mass error is smaller than EI, LE and LEEI2D methods.

\ExecuteMetaData[figures.tex]{fig:2drvcfls}

The change of geometrical error and mass error with the grid resolution are plotted in Fig. \ref{fig:2drvcfls}.
Again, when CFL number is small (CFL=0.05 as shown in Figs. \ref{fig:2drvcfls} (a) and (b)),
the geometrical error are similar to each other for all methods,
both WY and EILE2D method conserve mass to machine precision.
LEEI2D method show a better mass conservation than EI and LE methods.
When CFL number is large (CFL=1.0 as shown in Figs. \ref{fig:2drvcfls} (c) and (d)),
LEEI2D method has the best performance on geometrical error,
and better mass conservation than EI and LE method.
The EILE2D method ultimately conserve mass
and has the second best performance on geometrical error.
The WY method, which is similar performance on geometric
error,
no longer preserve mass when CFL = 1.
However, the WY has a better performance on mass conservation
than the LEEI2D method.

\subsection{3D tests}
\label{sec:3dtests}
\subsubsection{Zalesak's disk rotation test}
\label{sec:zalesak3d}

Similar to the 2D Zalesak's in subsection \ref{sec:zalesak2d},
the 3D Zalesak's roration test extends the slotted disk to a 3D
slotted sphere \citep{enright_hybrid_2002}.
The $x$ and $y$ component of the velocity remains the 2D rotation field
and a uniform velocity is applied to the $z$ component of the velocity field.
The initial setup and the velocity field are shown in Fig. \ref{fig:3dzalesakinit}.
Periodic boundary condition is applied at $z$ direction.

\ExecuteMetaData[figures.tex]{fig:3dzalesakinit}

\ExecuteMetaData[figures.tex]{fig:3dzalesakcfls}

\ExecuteMetaData[figures.tex]{fig:3dzalesakresos}

The results from different methods are plotted in
Fig. \ref{fig:3dzalesakcfls}.
Similar to the 2D result,
the geometrical error decreases when the CFL number becomes larger.
When the grid resolution increases, again, the geometrical error decreases.
Most of the methods do not show significant
difference between each other,
except for the EILE3D method.
This could be explained that
the velocity magnitude of three divergence-free 2D
is about a half of the original velocity field.
While the velocity magitude is halved,
the equavilent CFL numbers is halved as well,
leading to a larger geometrical error as the relationship
between the geomrtric error and CFL number shown in
Fig. \ref{fig:3dzalesakresos}.
For the mass error,
all results conserve to machine precision,
which remains the same as the 2D problem.

\subsubsection{Single vortex test}
\label{sec:sv3d}

The 2D Rider-Kothe single vortex problem is extended to 3D
with several approaches
\citep{liovic_3d_2006,aulisa_interface_2007,baraldi_mass-conserving_2014}.
The divergence-free velocity field of
3D single vortex test by \citet{liovic_3d_2006} is used in this study.
The initial set up and test parameters are shown in Fig. \ref{fig:3dsvinit}.
The $x$, $y$ component of the velocity field is identical
to the 2D Rider-Kothe single vortex test,
and the additional velocity component is a laminar pipe flow.

\ExecuteMetaData[figures.tex]{fig:3dsvinit}

The material interface at $t=T/2$ and $t=T$ with grid resolution $N=64$
is shown in Fig.
\ref{fig:3dsviso1} and Fig. \ref{fig:3dsviso2}.
The exact results in
Fig.\ref{fig:3dsviso1} (a) and Fig. \ref{fig:3dsviso2} (a)
are calculated with analytical marker particles.
At $t=T/2$,
all results does not have significant visual difference.
At $t=T$,
all results look similar with each other, except for the EI result.

\ExecuteMetaData[figures.tex]{fig:3dsviso1}
\ExecuteMetaData[figures.tex]{fig:3dsviso2}

\ExecuteMetaData[figures.tex]{fig:3dsvresos}

Fig. \ref{fig:3dsvresos} shows the effect on
the solution accuracy and mass conservation
with CFL number when the grid resolution is
$32 \times 32 \times 64$.
Using the velocity decomposition of
Eq. \eqref{eq:eile3ddecomposition}
on the velocity field in Fig. \ref{fig:3dsvinit},
the maximum magnitude of the velocity field
could be as big as about 3.0,
which means the equivalent maximum CFL is tripled
after velocity decomposition.
In Fig. \ref{fig:3dsvresos} (a),
when $\rm{CFL} \ge 0.5$,
unreasonable results lead to a spike on mass and geometrical error.
Among the rest of 5 methods,
EI method have most significant geometrical error,
while the geometrical error of LE methods is smaller than EI methods.
EILE3DS, EIEALE and WY methods has very similar performance on geometrical error
and better than those of EI and LE method.
For mass conservation,
EIEALE conserves mass with higher order of mass loss
than EILE3DS method.
WY method conserves mass to machine precision $\rm{CFL} < 0.5$,
for a larger CFL number,
the WY methods no longer conserves mass to machine precision.
However,
the mass error is very close to the EIEALE method and better than the rest of the methods.

\ExecuteMetaData[figures.tex]{fig:3dsvcfls}

The effect on
the solution accuracy and mass conservation
with repeat halving the grid resolution
is shown in Fig. \ref{fig:3dsvcfls}.
Note that most of the EILE3D results are out of the plot
range.
With different grid resolution under two CFL numbers (0.3 and 1.0).
Considering that the velocity decomposition leads to
the stability issues on EILE3DS methods,
WY method shows an overall best performance on
mass conservation.

\subsubsection{Reverse Vortex test}
\label{sec:rv3d}

A more stringent 3D test: reversed vortex test
\citep{leveque_high-resolution_1996} is taken in this subsection.
The initial set up and test parameters are shown in Fig. \ref{fig:3drvinit}.
Under the two rotating vortices,
the initial material of sphere scoops out and pancake,
leading to a very thin and stretched interface at $t=T/2$
(See Fig. \ref{fig:3drviso1}(a)),
and recovers to the initial sphere at $t=T$
(See Fig. \ref{fig:3drviso2}(a)).
When the grid resolution is not fine enough,
part of the interface may thin out of the grid resolution,
leading to a massive of fragmentation and coalescence.
When recovered the initial position at $t=T$
(See Fig. \ref{fig:3drviso1} and Fig. \ref{fig:3drviso2}).

\ExecuteMetaData[figures.tex]{fig:3drvinit}

\ExecuteMetaData[figures.tex]{fig:3drviso1}
\ExecuteMetaData[figures.tex]{fig:3drviso2}

\ExecuteMetaData[figures.tex]{fig:3drvresos}

Fig. \ref{fig:3drvresos} shows the effect on
the solution accuracy and mass conservation
with CFL number when the grid resolution is
$64 \times 64 \times 64$.
Unlike the velocity field in Section \ref{sec:sv3d},
when applying the velocity decomposition of
Eq. \eqref{eq:eile3ddecomposition}
applied on the Fig. \ref{fig:3dsvinit},
the maximum magnitude of the velocity
is about 0.5,
which is about a half of the original velocity field.
The velocity decomposition of Eq. \eqref{eq:eile3ddecomposition}
will not bring in the
stability issue for the velocity field shown in
Fig. \ref{fig:3drvinit}.
The EILE3D methods always conserve mass to machine precision
and performs as one of the best methods among all six
candidate methods on geometrical error.
The EIEALE and WY methods has similar performance on
geometrical error.
When $\rm{CFL} \le 0.5$,
WY method conserves mass to machine precision,
and has a slightly smaller mass error compared with the
EILE method.
The EI and LE methods again has the worst mass conservation,
and LE method has smaller geometrical error than the EI methods.
the EILE3D method,
although improves the mass the geometrical error over EI
and LE methods,
is not as good as EIEALE, WY or EILE3D method in this test.

\ExecuteMetaData[figures.tex]{fig:3drvcfls}

The effect on
the solution accuracy and mass conservation
with repeat halving the grid resolution
is shown in Fig. \ref{fig:3drvcfls}.
The EILE3D method has the overall best performance on mass conservation and geometrical error.
The WY method has similar performance with the EILE3D method,
except for the mass conservation on large CFL number.
However, even when unconserved mass occur due to the
large CFL number,
WY method is one of the best method on mass conservation among the rest of the methods exclude EILE3D method.
WY method improves the mass conservation by the order of
$O(10^{-2}) - O(10^{-4})$ compared with EI and LE methods.

\section{Conclusion}
\label{sec:conclusion}

In this study,
a comparative study on several flux-splitting MOF methods are carried out.
Several directional splitting methods has been extended from VOF-PLIC to MOF
advection.
These methods, along with several existing split MOF advection methods,
are investigated in both 2D and 3D.
We focus on the three important conditions for MOF advection:
(1) global mass conservation condition,
(2) local bounded condition
(3) geometrical consistency condition.
Six numerical tests cases are used to evaluate the geometrical and mass
conservation
for different MOF advection methods.
Besides, we also derived a fast and efficient formulation for calculating the volume and
centroid
from known equation of plane.

For 2D problem,
the hybrid Eulerian Implicit and Lagrangian Explicit (EI-LE) strategy
is the only strategy that satisfies all three conditions.
The EI-LE strategy in 2D (EILE2D) is ultimately mass conservation and geometrical
conservation.
The EILE2D method has also obtained the overall best performance on
geometrical error.
For 2D problem,
EILE2D method is recommended for MOF advection.

The implementation of EI-LE strategy in 3D (EILE3D) is not as straightforward as
its in 2D.
The EILE3D method need to decomposed the original velocity field into 3 2D
divergence-free velocity fields,
which is equivalent with 3 EILE2D advection.
Thus three conditions holds for EILE3D.
However,
the numerical tests show that the velocity decomposition may
bring in stability issue.
The Weymouth-Yue strategy (WY) is ultimately mass conservation,
although it does not guarantee the geometrical consistency,
the overall performance on geometrical error is one of the best among all split
advection methods in this study.
While the WY method has a more stringent CFL limit of 1/6 in order to guarantee
an ultimate mass conservation.
Numerical tests show that for a larger CFL number,
although WY method no longer preserve mass to machine precision,
it still significantly improves the mass conservation from EI or LE methods.
Considering the overall performance on mass conservation and geometrical error,
along with the stability issue,
WY methods is preferable for 3D problems.




\bibliography{MOF-EILE}

\newpage

\section*{Appendix A: Triangular cutting plane}
\subsection*{Configuration 1: Triangle cut interface}
This configuration is the simplest one as all intercepts of the $x_i$ coordinates
are within the range of the hexahedron.
The cutting interface is a triangle.
The volume and corresponding centroid are expressed as
\begin{equation}
    \label{eq:cut1V}
    V=\frac{\alpha^{3}}{6 m_{1} m_{2} m_{3}}
\end{equation}
and the corresponding centroid is expressed as
\begin{equation}
    \begin{aligned}
        \label{eq:cut1C}
        c_1=\frac{\alpha}{4 m_{1}}, \\
        c_2=\frac{\alpha}{4 m_{2}}, \\
        c_3=\frac{\alpha}{4 m_{3}}.
    \end{aligned}
\end{equation}

\subsection*{Configuration 2: Quadrilateral cutting interface}
When $m_1 < \alpha < m_3$, the intercept on $x_1$ coordinate exceeds the range of
of the hexahedron.
The cutting interface is a quadrilateral.
The cutting polyhedron can be determined by subtracting a smaller tetrahedron from the
tetrahedron formed by the original point and three intercepts points.
The volume and corresponding centroid are expressed as
\begin{equation}
    \label{eq:cut2V}
    V=\frac{A}
    {6 m_{1} m_{2} m_{3}},
\end{equation}
\begin{equation}
    \begin{aligned}
        \label{eq:cut2C}
        c_1 & =\frac{3}{4} - \frac{\alpha\left(3 \alpha-m_{1}\right)}{4 A},          \\
        c_2 & =\frac{\left(2 \alpha-m_{1}\right)(\alpha^2 + (\alpha-m_{1})^2)}{4 A}, \\
        c_3 & =\frac{\left(2 \alpha-m_{1}\right)(\alpha^2 + (\alpha-m_{1})^2)}{4 A},
    \end{aligned}
\end{equation}
where
\begin{equation}
    \begin{aligned}
        \label{eq:cut2ext}
        A = & \alpha^3 - (\alpha - m_1)^3 - (\alpha - m_2)^3,
    \end{aligned}
\end{equation}

\subsection*{Configuration 3: Pentagon cutting interface}
When $m_2 < \alpha < m_c$, where $m_c = \min(m_{12},m_3)$ and $m_{12} = m_1 + m_2$,
the intercept on $x_1$ and $x_2$ coordinates exceeds the range of of the hexahedron.
The cutting interface is a pentagon.
The cutting polyhedron can be determined by subtracting two smaller tetrahedron from the
tetrahedron formed by the original point and three intercepts points.
The volume and corresponding centroid are expressed as
\begin{equation}
    \label{eq:cut3V}
    V =\frac{A}{6 m_{1} m_{2} m_{3}},
\end{equation}
\begin{equation}
    \begin{aligned}
        \label{eq:cut3C}
        c_1 & =-\frac{\alpha^4 - (\alpha-m_1)^3(\alpha+3m_1) - (\alpha-m_2)^4}{4m_1A}, \\
        c_2 & =-\frac{\alpha^4 - (\alpha-m_1)^4 - (\alpha-m_2)^3(\alpha+3m_2)}{4m_2A}, \\
        c_3 & =-\frac{\alpha^4 - (\alpha-m_1)^4 - (\alpha-m_2)^4}{4m_3A},              \\
    \end{aligned}
\end{equation}
where
\begin{equation}
    \begin{aligned}
        \label{eq:cut3ext}
        A = & \alpha^3 - (\alpha-m_1)^3 - (\alpha - m_2)^3.
    \end{aligned}
\end{equation}

\subsection*{Configuration 4: Hexagon cutting interface}
When $m_3<m_{12}$ and $m_3 < \alpha < 1/2$,
all intercepts on $x_i$ coordinates exceeds the range of of the hexahedron.
The cutting interface is a hexagon.
The cutting polyhedron can be determined by subtracting three asmaller tetrahedron from the
tetrahedron formed by the original ponit and three intercepts points.
The volume and corresponding centroid are expressed as
\begin{equation}
    \label{eq:cut4V}
    \mathrm{V}=\frac{A}{6 m_{1} m_{2} m_{3}},
\end{equation}

\begin{equation}
    \begin{aligned}
        \label{eq:cut4C}
        c_1 & =-\frac{\alpha^4 - (\alpha-m_1)^3(\alpha+3m_1) - (\alpha-m_2)^4 - (\alpha-m_3)^4}{4m_1A}, \\
        c_2 & =-\frac{\alpha^4 - (\alpha-m_1)^4 - (\alpha-m_2)^3(\alpha+3m_2) - (\alpha-m_3)^4}{4m_2A}, \\
        c_3 & =-\frac{\alpha^4 - (\alpha-m_1)^4 - (\alpha-m_2)^4 - (\alpha-m_3)^3(\alpha+3m_3)}{4m_3A}, \\
    \end{aligned}
\end{equation}
where
\begin{equation}
    \begin{aligned}
        \label{eq:cut4ext}
        A = & \alpha^3 - (\alpha-m_1)^3 - (\alpha - m_2)^3 - (\alpha - m_3)^3. \\
    \end{aligned}
\end{equation}

\subsection*{Configuration 5: Quadrilateral cutting interface 2}
When $m_{12}<m_{3}$ and $m_{12} < \alpha < 1/2$,
all intercepts on $x_i$ coordinates exceeds the range of of the hexahedron.
The cutting interface is a hexagon.
The cutting polyhedron can be determined by subtracting two smaller tetrahedron from the
tetrahedron formed by the original ponit and three intercepts points,
then combine the intersection of the two small tetrahedron.
The volume and corresponding centroid are expressed as
\begin{equation}
    \label{eq:cut5V}
    V =\frac{A}{6 m_{1} m_{2} m_{3}},
\end{equation}

\begin{equation}
    \begin{aligned}
        \label{eq:cut5C}
        c_1 & = \frac{1}{2} - \frac{m_1}{6(2\alpha-m_{12})},                                      \\
        c_2 & = \frac{1}{2} - \frac{m_2}{6(2\alpha-m_{12})},                                      \\
        c_3 & = \frac{3(\alpha-2m_{12})(2\alpha-m_{12})+\alpha m_{12}-m_1m_2}{6(2\alpha-m_{12})}, \\
    \end{aligned}
\end{equation}
where
\begin{equation}
    \begin{aligned}
        \label{eq:cut5ext}
        A = & 3 m_1 m_2 (2\alpha - m_1 - m_2).
    \end{aligned}
\end{equation}

\newpage

\section*{Appendix B: Tables for numerical tests}

\ExecuteMetaData[tables.tex]{tab:2dzalesakcfls}

\ExecuteMetaData[tables.tex]{tab:2dzalesakresos}

\ExecuteMetaData[tables.tex]{tab:2dsvresos}

\ExecuteMetaData[tables.tex]{tab:2dsvresosmass}

\ExecuteMetaData[tables.tex]{tab:2dsvcfls}

\ExecuteMetaData[tables.tex]{tab:2dsvcflsmass}

\ExecuteMetaData[tables.tex]{tab:2drvresos}

\ExecuteMetaData[tables.tex]{tab:2drvresosmass}

\ExecuteMetaData[tables.tex]{tab:2drvcfls}

\ExecuteMetaData[tables.tex]{tab:2drvcflsmass}

\ExecuteMetaData[tables.tex]{tab:3dzalesakcfls}

\ExecuteMetaData[tables.tex]{tab:3dzalesakresos}

\ExecuteMetaData[tables.tex]{tab:3dsvresos}

\ExecuteMetaData[tables.tex]{tab:3dsvresosmass}

\ExecuteMetaData[tables.tex]{tab:3dsvcfls}

\ExecuteMetaData[tables.tex]{tab:3dsvcflsmass}

\ExecuteMetaData[tables.tex]{tab:3drvresos}

\ExecuteMetaData[tables.tex]{tab:3drvresosmass}

\ExecuteMetaData[tables.tex]{tab:3drvcfls}

\ExecuteMetaData[tables.tex]{tab:3drvcflsmass}

\end{document}